\documentclass[graybox]{svmult}

\usepackage{mathptmx}       

\usepackage[utf8]{inputenc}

\usepackage{helvet}         
\usepackage{courier}        
\usepackage{type1cm}        
\usepackage{makeidx}         
\usepackage{graphicx}        
\usepackage{multicol}        
\usepackage[bottom]{footmisc}

\makeindex             

\usepackage{csquotes}
\usepackage{afterpage}
\usepackage[T1]{fontenc}
\usepackage{lmodern}
\usepackage{amscd,indentfirst,epsfig}
\usepackage{times}
\usepackage{enumerate}
\usepackage{mathrsfs}
\usepackage{stmaryrd}
\usepackage{amsopn}
\usepackage{empheq,bm}
\usepackage{amsbsy,amsmath}
\usepackage{amscd}
\usepackage{comment}
\usepackage{appendix}
\usepackage{caption}
\usepackage{microtype}
\usepackage{mathtools, amsthm, amsfonts, amssymb, stmaryrd, dsfont}
\usepackage{caption}
\usepackage{subcaption}
\usepackage{url}
\usepackage{doi}
\usepackage{csquotes}
\usepackage{colortbl}
\usepackage{ esint }
\usepackage{arydshln}
\usepackage{tabularx}
\usepackage{chngcntr}
\usepackage{etoolbox}
\usepackage{array}
\usepackage{tikz}
\usepackage[most]{tcolorbox}
\usepackage{xcolor}
\usepackage{listings}
\usetikzlibrary{shadows}
\usepackage{subcaption}

\usepackage{bm}

\allowdisplaybreaks

\usepackage{mathptmx}    

\tcbset{colback=yellow!10!white, colframe=red!50!black, 
        highlight math style= {enhanced, 
            colframe=red,colback=red!10!white,boxsep=0pt}
        }

\newtcolorbox{mybox}[1]{colback=blue!5!white,colframe=blue!50,fonttitle=\bfseries,title=#1}

\renewcommand{\arraystretch}{1.5}

\definecolor{myblue}{rgb}{.9, .9, 1}

\usetikzlibrary{arrows}
\usetikzlibrary{shapes.geometric,decorations.markings,decorations.pathmorphing,decorations.pathreplacing,calligraphy}
\usetikzlibrary{patterns}
\usetikzlibrary{shadows}

\newcommand{\AxisRotator}[1][rotate=0]{%
    \tikz [x=0.20cm,y=0.40cm,line width=.2ex,-stealth,#1] \draw (0,0) arc (-150:150:1 and 1);%
}

\tikzset{
diagonal fill/.style 2 args={fill=#2, path picture={
\fill[#1, sharp corners] (path picture bounding box.south west) -|
                         (path picture bounding box.north east) -- cycle;}},
reversed diagonal fill/.style 2 args={fill=#2, path picture={
\fill[#1, sharp corners] (path picture bounding box.north west) |- 
                         (path picture bounding box.south east) -- cycle;}}
}

\allowdisplaybreaks

\newcommand{\R}{\mathbb{R}}
\newcommand{\N}{\mathbb{N}}

\newcommand{\J}{\mathcal{J}}
\renewcommand{\E}{\mathcal{E}}

\newcommand{\dd}{\mathrm{d}}

\newcommand{\e}{\varepsilon}

\renewcommand{\leq}{\leqslant}
\renewcommand{\geq}{\geqslant}

\begin{document}

\title*{On the modelling of polyatomic molecules in kinetic theory}
\author{M. Bisi, T. Borsoni \& M. Groppi}
\institute{M. Bisi\at Università di Parma, Parma, Italy. 
\email{marzia.bisi@unipr.it}
\and T. Borsoni \at CERMICS, École Nationale des Ponts et Chaussées, Champs-sur-Marne, France.
\email{thomas.borsoni@enpc.fr}
\and
M. Groppi$^*$ \at Università di Parma, Parma, Italy. 
\email{maria.groppi@unipr.it}}

%
%
\maketitle

\abstract{This communication is both a pedagogical note for understanding polyatomic modelling in kinetic theory and a ``cheat sheet'' for a series of corresponding concepts and formulas. We explain, detail and relate three possible approaches for modelling the polyatomic internal structure, that are: the internal states approach, well suited for physical modelling and general proofs, the internal energy levels approach, useful for analytic studies and corresponding to the common models of the literature, and the internal energy quantiles approach, less known while being a powerful tool for particle-based numerical simulations such as Direct Simulation Monte-Carlo (DSMC). 
This note may in particular be useful in the study of non-polytropic gases.}



\section*{Introduction}

The construction of a statistical model for a system of molecules, such as Boltzmann, BGK or Fokker-Planck, requires both a \textbf{single-molecule} and an \textbf{inter-molecular interaction} description.

\medskip

\noindent The \textbf{single-molecule} model determines the \emph{variables} which the studied molecular density should depend on, as well as the \emph{structure of the space} they belong to -- typically, the \emph{measure} of which their belonging space is endowed with. On the other hand, the \textbf{inter-molecular interaction} model shall determine the \emph{form} of the operators involved in the evolution equation of the molecular density.

\medskip

\noindent A precise single-molecule model is therefore an unavoidable first step for a precise statistical model. The purpose of this short note is to provide a pedagogical summary of the necessary tools (insights and formulas) for the mathematical setting of a \textbf{polyatomic} single-molecule model  designed for kinetic theory, related to various previous works of the community (\cite{arima2018extended, bisi2005kinetic, borgnakke1975statistical,desvillettes1997modele, gamba2020cauchy, pavic2022kinetic, wang1951transport} to cite a few), and based on works of the authors~\cite[Introduction, Section 1.2]{borsoni2024contributions} and~\cite{GF}.
\medskip

\noindent Non-relativistic polyatomic gases are composed by particles characterized by their own velocity and by an additional internal state energy variable, representing the non-translational degrees of freedom, that may be assumed to be discrete or continuous.  Wang Chang and Uhlenbeck \cite{wang1951transport} in the Fifties  were the first  to describe the internal structure of a polyatomic gas using a discrete energy variable, related to the vibrational modes. This approach has been extensively studied in kinetic theory around the  2000s by Giovangigli~\cite{giovangigliBook}, Spiga, Groppi and Bisi \cite{bisi2005kinetic,groppi1999kinetic}. On the other hand, Borgnakke and Larsen \cite{borgnakke1975statistical} in the Seventies  proposed a description of the internal energy structure based on a continuous variable, suitable to represent rotational modes. This approach has been then investigated by Desvillettes \cite{desvillettes1997modele}, showing that the internal structure is in fact specified by the choice of the measure associated to the energy variable. Suitable choices of such a measure allow the description of non-polytropic gases, for which the internal energy depends on the temperature in a non-linear way.
\medskip

\noindent In this note we discuss and compare three possible general models  for the polyatomic internal structure. The first is the internal states approach, which is well suited for physical modelling and general proofs, investigated in detail in \cite{GF}. The second is the internal energy levels approach, useful for analytic studies and corresponding to the common models of the literature. The third is the internal energy quantiles approach, which is a simpler and  powerful tool for particle-based numerical simulations such as Direct Simulation Monte-Carlo (DSMC).

\medskip

\noindent In Section~\ref{section1}, we focus on the formulation of the  polyatomic single-molecule model, then we detail an example for a diatomic molecule in Section~\ref{section2} and discuss in Section~\ref{section3} how to incorporate the single-molecule model at the kinetic level.

\section{Modelling a polyatomic molecule} \label{section1}

The field focused on molecular models is known as Molecular Modelling (see for instance~\cite{leach2001molecular}) and has a wide range of applications. While a comprehensive model of a specific molecule is not in the scope of this note, our objective here is to set a mathematical framework that is both physically accurate, adaptable to molecular models and mathematically manageable in kinetic theory.

\medskip

\subsection{A general description based on internal states}

\smallskip

We consider the molecule to have, at each time instant $t$, a position, a velocity, and an \emph{internal state}, respectively called $x$, $v$ and \textcolor{orange}{$\zeta$}. In the following, we focus on \textcolor{orange}{$\zeta$} and discard the variables $x$ and $v$.

\smallskip

\begin{minipage}{.23\textwidth}
  \hspace{-2em}  \begin{tikzpicture}[scale=.23,decoration={coil,aspect=0.5,segment length = 4mm, amplitude=3mm}]
\tikzset{molecule_normal/.pic={code= \tikzset{scale=.25}
    \shade[ball color=blue,opacity=1] (-1.6,3.2) circle (1.2);
    \shade[ball color=red,opacity=1] (0,4) circle (2);
    \shade[ball color=blue,opacity=1] (1.6,2.7) circle (1.2);
    }}
\tikzset{molecule/.pic={code= \tikzset{scale=.3}
    \shade[ball color=black,opacity=1] (-1.6,3.2) circle (1.2);
    \shade[ball color=yellow!20,opacity=1] (0,4) circle (2);
    \shade[ball color=black,opacity=1] (1.6,2.7) circle (1.2);
    }}

\tikzset{ressort/.style={thin,orange!60!gray,smooth}}

\node[circle, inner color = orange, outer color = white]
    at (-15,10) {$\hspace{45pt}$};
\node[circle, inner color = orange, outer color = white]
    at (-20,2) {$\hspace{25pt}$};
    \node[circle, inner color = orange, outer color = white]
    at (-11,0) {$\hspace{25pt}$};
\shade[ball color=red] (-15,10) circle (2);
\shade[ball color=blue] (-20,2) circle (1.2);
\shade[ball color=blue] (-11,0) circle (1.2);
\draw[very thick, opacity=.9] (-19.7,2.5)--(-16,8.3) ;
\draw[very thick, opacity=.9] (-11.4,1.05)--(-14.5,8.8) ;

\draw[ressort, decorate] (-19.7,2.5)--(-16,8.3) ;
\draw[ressort, decorate] (-11.4,1.05)--(-14.5,8.8) ;
\draw[ressort, decorate] (-19.5,1.8)--(-12.18,.2) ;

\draw[orange] [rotate around={160:(-14.7,3.6)}, line width=1.5pt, -stealth, dashed] (-10,2) arc (-30:210:4.5cm and 1cm);

\node at (-11,8) {\textcolor{orange}{$\zeta$}};

\end{tikzpicture}
\end{minipage}
\begin{minipage}{.72\textwidth}
    
 The internal state variable \textcolor{orange}{$\zeta$} would typically encompass rotation and vibrations, nevertheless it may also contain other information such as orientation, electronic structure... It is assumed to  belong to a \emph{space of internal states} $\E$, which is equipped with a measure $\mu$, the \emph{degeneracy} of the states. 

 \medskip
 
 To each internal state $\zeta$ corresponds an \emph{internal energy}, denoted by $\e(\zeta)$. Any molecule has a \emph{fundamental energy level}, called $\e^0$, which is the minimum of the function $\e$.

\end{minipage}

\medskip

\begin{example}
\emph{Consider a non-linear triatomic molecule for which we model separately rotation via classical mechanics (rigid rotor) and vibrations via the quantum harmonic oscillator approach. The internal state $\zeta$ of the molecule is in this case composed of the angular velocity $\omega \in \R^3$ (in a reference frame keeping the molecule's orientation fixed) and vibration modes $(n_1,n_2,n_3) \in \N^3$. The space of internal states is therefore $\E = \R^3 \times \N^3$, equipped with $\mu = \mathrm{Lebesgue}_{\R^3} \otimes \mathrm{Counting}_{\N^3}$, 
 avoiding degeneracy.}
\end{example}

\hspace{-1.4em} \begin{minipage}{.79\textwidth}
\noindent \emph{A molecule in the state $\zeta = (\omega,n)$ bears an energy $\e(\zeta)$, equal to the sum of its kinetic energy of rotation and energy associated with each mode of vibration, that is}
\[
\e(\omega,n) = \frac12 \; \sum_{i=1}^3 \mathcal{J}_i \, \omega_i^2 + \sum_{i=1}^3 \Delta \epsilon_i \; \left(n_i + \frac12 \right),
\]
\emph{with $\mathcal{J} \in \R_+^3$ the moment of inertia and $\Delta \epsilon_i$ the gap between two energy levels of the $i^{th}$ vibration mode. The fundamental energy level is $\e^0 = \frac12 (\Delta \epsilon_1 + \Delta \epsilon_2 + \Delta \epsilon_3)$.}
\end{minipage}
\begin{minipage}{.2\textwidth}
   \hspace{-1em} \begin{tikzpicture}[scale=.2,decoration={coil,aspect=0.5,segment length = 4mm, amplitude=3mm}]
\tikzset{molecule_normal/.pic={code= \tikzset{scale=.25}
    \shade[ball color=blue,opacity=1] (-1.6,3.2) circle (1.2);
    \shade[ball color=red,opacity=1] (0,4) circle (2);
    \shade[ball color=blue,opacity=1] (1.6,2.7) circle (1.2);
    }}
\tikzset{molecule/.pic={code= \tikzset{scale=.3}
    \shade[ball color=black,opacity=1] (-1.6,3.2) circle (1.2);
    \shade[ball color=yellow!20,opacity=1] (0,4) circle (2);
    \shade[ball color=black,opacity=1] (1.6,2.7) circle (1.2);
    }}

\tikzset{ressort/.style={thin,orange!40!gray,smooth}}
\shade[ball color=red] (-15,10) circle (2);
\shade[ball color=blue] (-20,2) circle (1.2);
\shade[ball color=blue] (-11,0) circle (1.2);
\draw[very thick, opacity=.9] (-19.7,2.5)--(-16,8.3) ;
\draw[very thick, opacity=.9] (-11.4,1.05)--(-14.5,8.8) ;

\draw[ressort, decorate] (-19.7,2.5)--(-16,8.3) ;
\draw[orange!40!gray] (-21,5) node{$n_1$};
\draw[ressort, decorate] (-11.4,1.05)--(-14.5,8.8) ;
\draw[orange!40!gray] (-11,6.5) node{$n_2$};
\draw[ressort, decorate] (-19.5,1.8)--(-12.18,.2) ;
\draw[orange!40!gray] (-16,-1) node{$n_3$};

\draw[orange] [rotate around={160:(-14.7,3.6)}, line width=1.5pt, -stealth, dashed] (-10,2) arc (-30:210:4.5cm and 1cm);
\draw[orange] (-20,9) node{$\omega$};

\end{tikzpicture}
\end{minipage}

\medskip

\medskip

\noindent In the above example, other models could have been considered to describe rotation and vibration, or choosing to consider solely rotation, or taking subtler effects into account...

\begin{mybox}{To sum up}

\noindent An \textbf{internal states description} of a molecule is constructed with
\begin{align*}
\text{-- } &\text{a space of internal states } \E \text{ equipped with a measure } \mu \\
\text{-- } &\text{an internal energy function } \e : \E \to \R \text{ measurable, which minimum is } \e^0.
\end{align*}
The internal variable is the \emph{internal state} $\zeta \in \E$, its associated energy is $\e(\zeta)$.
\end{mybox}

\subsection{Descriptions based on internal energy}

While the \textbf{internal states} setting presented above is fundamental, it is relevant to consider a setting where the information is compacted, based on \textbf{internal energy}. We  present the two main descriptions. One is  based on internal energy \emph{levels} and the other on internal energy \emph{quantiles}.

\subsubsection{Energy \emph{levels}} \label{section:energy_levels}

First of all, we highlight that the energy levels approach presented hereafter is in correspondence with the common models of the literature, the model with continuous energy~\cite{borgnakke1975statistical,desvillettes1997modele}  and the model with discrete energy~\cite{bisi2005kinetic,wang1951transport}, and this is explained later on.

\medskip

Let us build an ``internal energy levels'' model from an ``internal states'' one, which space of states is $\E$ endowed with $\mu$, and internal energy function is $\e$. We let\footnote{At the kinetic level, the value of $\e^0$ matters only when chemical reactions are involved. For anything that is not directly related to chemical reactions, it is handier to consider $\bar{\e}$ rather than $\e$.}
\begin{equation}\bar{\e} := \e - \e^0\end{equation} to be the \emph{grounded} internal energy function, which image is contained in $\R_+$ and essential infimum is $0$.  Define the \textbf{energy law} as the image measure, on $\R_+$, of $\mu$ by $\Bar{\e}$,
\begin{equation} \label{eqdef:mesureimagemubareps_intro}
   \qquad \qquad   \mu_{\Bar{\e}}(B) := \mu \left(\Bar{\e} \in B \right), \qquad B \in \mathrm{Bor}(\R_+).
\end{equation}
Then, for any test function $\phi$,
\[
\int_{\E} \phi \big(\varepsilon  (\zeta) \big) \, \dd\mu(\zeta)  = \int_{\R_+} \phi (I + \varepsilon^0) \, \dd\mu_{\Bar{\varepsilon}}(I).
\]
Identifying, there is therefore a correspondence\footnote{We see in a further section in which situation we can consider these descriptions to be equivalent.} between the base \emph{internal state} model and a model in which
\begin{enumerate}
    \item we consider the variable $I$, called \textbf{energy level}, in the space $\R_+$ equipped with $\mu_{\Bar{\e}}$,
    \item the internal energy of the molecule is $I + \e^0$.
\end{enumerate}
This is called an \textbf{internal energy levels} description.

\medskip

\noindent Notice that if $\mu_{\Bar{\e}}$ has a density with respect to the Lebesgue measure, \textbf{which is called} $\varphi$ \textbf{in the literature}, we recover the model with continuous energy~\cite{borgnakke1975statistical,desvillettes1997modele} (taking $\e_0 = 0$). \textbf{We highlight that~\eqref{eqdef:mesureimagemubareps_intro} indicates how to compute $\varphi$ from the molecular model}. In the case where $\mu_{\Bar{\e}}$ is a discrete measure, the discrete energy model~\cite{bisi2005kinetic,wang1951transport} is recovered.

\begin{mybox}{To sum up}

\noindent An \textbf{internal energy levels description} of a molecule is constructed with
\begin{align*}
\text{-- } &\text{an energy law } \mu_{\bar{\e}} \text{ on } \R_+ \\
\text{-- } &\text{a fundamental internal energy level } \e^0.
\end{align*}
The internal variable is the \emph{internal energy level} $I \in \R_+$, its associated energy is $I + \e^0$.
\end{mybox}

\subsubsection{Energy \emph{quantiles}} \label{section:energy_quantiles}

Let us build an ``energy quantiles'' model from an ``energy levels'' one, with an energy law $\mu_{\bar{\e}}$ and a fundamental energy level $\e^0$. Define the cumulative internal energy function $F_{\mu_{\Bar{\varepsilon}}}$ on $\R_+$ by
\begin{equation} \label{eqdef:fonctionderepartition_intro}
    F_{\mu_{\Bar{\varepsilon}}}(I) :=  \mu_{\Bar{\varepsilon}} \left([0, \,I) \right), \qquad I \in \R_+,
\end{equation}
and the \emph{energy quantile function} $F_{\mu_{\Bar{\varepsilon}}}^{\leftarrow}$ on $(0,q^{\max})$ as the left generalised inverse of $F_{\mu_{\Bar{\e}}}$, by
\begin{equation} \label{eqdef:fonctionquantile_intro}
    F_{\mu_{\Bar{\varepsilon}}}^{\leftarrow}(q) := \inf \left\{ I \in \R_+ \text{ such that } F_{\mu_{\Bar{\varepsilon}}}(I) \geq q  \right\}, \qquad q \in ( 0, q^{\max}),
\end{equation}
where $q^{\max} = \mu_{\Bar{\varepsilon}}(\R_+)$. Of course, when $F_{\mu_{\Bar{\varepsilon}}}$ is invertible, $F_{\mu_{\Bar{\varepsilon}}}^{\leftarrow} = F_{\mu_{\Bar{\varepsilon}}}^{-1}$.

\medskip

\noindent Then, for any test function $\phi$,
$$
\int_{\R_+} \phi (I + \e^0 ) \, \dd\mu_{\Bar{\varepsilon}}(I)  = \int_{(0,q^{\max})} \phi \left(F_{\mu_{\Bar{\varepsilon}}}^{\leftarrow}(q) + \e^0 \right) \, \dd q.
$$

\medskip

\noindent Identifying, there is therefore a correspondence\footnote{The two approaches are always equivalent.} between the base \emph{internal energy level} model and a model in which
\begin{enumerate}
    \item we consider the variable $q$, called \textbf{energy quantile}, in the space $(0,q^{\max})$, equipped with the \emph{Lebesgue} measure,
    \item the internal energy of the molecule is $F_{\mu_{\Bar{\varepsilon}}}^{\leftarrow}(q) + \e^0$.
\end{enumerate}
This is called an \textbf{internal energy quantiles} description.

\begin{mybox}{To sum up}

\noindent An \textbf{internal energy quantiles description} of a molecule is constructed with
\begin{align*}
&\text{-- a quantiles interval } (0,q^{\max}), \text{ equipped with the } \mathrm{Lebesgue} \text{ measure}, \\
&\text{-- an internal energy quantile function } F^{\leftarrow}_{\mu_{\bar{\e}}} : (0,q^{\max}) \to \R_+  
\text{ and a} 
\\
&\text{\phantom{--} fundamental internal energy level } \e^0.
\end{align*}
The internal variable is the \emph{internal energy quantile} $q \in (0,q^{\max})$, its associated energy is $F^{\leftarrow}_{\mu_{\bar{\e}}}(q) + \e^0$.
\end{mybox}

\subsection{Links between descriptions} \label{section:proba}

The polyatomic setting presented above has a structure of a probability setting (see Table \ref{table:parallel}), the sole difference being that $\mu(\E)$ is not constrained to be $1$ (in particular, it can be infinite).

\begin{table}[h!] 
    \centering
  \def\arraystretch{1.8}
  \begin{tabularx}{\textwidth}{ 
  |>{\hsize=1\hsize\linewidth=\hsize\centering}X 
  | >{\hsize=1\hsize\centering\arraybackslash}X |}
  \hline
\cellcolor{black!15!white} \textbf{Polyatomic description} & \cellcolor{black!15!white}\textbf{Probability setting} \\
     \hline 
\cellcolor{orange!20!white}  $(\E,\mu)$ space of \textcolor{orange!75!black}{internal states}  & \cellcolor{orange!20!white}  $(\Omega, \mathbb{P})$ space of \textcolor{orange!75!black}{events} \\
  \cellcolor{orange!20!white}    $\bar{\e} \, $ (grounded) internal energy function & \cellcolor{orange!20!white}  $X$ real random variable \\
     \hline
\cellcolor{olive!10!white} $(\R_+, \, \mu_{\Bar{\varepsilon}})$ space of \textcolor{olive!75!black}{internal energy levels}  & \cellcolor{olive!10!white} $(\R, \, \mathbb{P}_X)$ space of \textcolor{olive!75!black}{realisations}   \\
   \cellcolor{olive!10!white}  $\mu_{\Bar{\varepsilon}}$ internal energy law & \cellcolor{olive!10!white} $\mathbb{P}_X$ law of $X$ \\
     \hline
 \cellcolor{gray!10!white}    $\left(( 0, \, q^{\max} ), \, \mathrm{Lebesgue} \right)$  space of  \textcolor{gray!70!black}{internal energy quantiles} &  \cellcolor{gray!10!white}   $\left((0,1), \,  \mathrm{Lebesgue} \right)$ space of  \textcolor{gray!70!black}{quantiles} \\
  \cellcolor{gray!10!white}     $F_{\mu_{\Bar{\varepsilon}}}^{\leftarrow}$ internal energy quantile function &  \cellcolor{gray!10!white}   $F^{\leftarrow}_{\mathbb{P}_X}$ quantile function \\
    \hline
\end{tabularx}
    \caption{Correspondence  between the internal framework and a probability setting.}
    \label{table:parallel}
\end{table}

With this interpretation at hand, the links between the descriptions are all the more straightforward:
\begin{align}
    \mu_{\bar{\e}} = \; &\bar{\e} \, \# \, \mu, \\
    \mu(\E) = \mu_{\bar{\e}}&(\R_+) = q^{\max}, \\
    F^{\leftarrow}_{\mu_{\bar{\e}}} \text{ left generalized } &\text{inverse of } I \mapsto \mu_{\bar{\e}}([0,I)).
\end{align}

\subsection{Summary and use of the three descriptions}

\begin{figure}[htbp]
    \centering

\begin{tikzpicture}[scale=.3,decoration={coil,aspect=0.5,segment length = 4mm, amplitude=3mm}]
\tikzset{molecule_normal/.pic={code= \tikzset{scale=.3}
    \shade[ball color=blue,opacity=1] (-1.6,3.2) circle (1.2);
    \shade[ball color=red,opacity=1] (0,4) circle (2);
    \shade[ball color=blue,opacity=1] (1.6,2.7) circle (1.2);
    }}
\tikzset{molecule/.pic={code= \tikzset{scale=.3}
    \shade[ball color=black,opacity=1] (-1.6,3.2) circle (1.2);
    \shade[ball color=yellow!20,opacity=1] (0,4) circle (2);
    \shade[ball color=black,opacity=1] (1.6,2.7) circle (1.2);
    }}
\tikzset{molecule_full/.pic={code= \tikzset{scale=.3}

\tikzset{ressort/.style={thin,orange!40!gray,smooth}}
\shade[ball color=red] (-15,10) circle (2);
\shade[ball color=blue] (-20,2) circle (1.2);
\shade[ball color=blue] (-11,0) circle (1.2);
\draw[very thick, opacity=.9] (-19.7,2.5)--(-16,8.3) ;
\draw[very thick, opacity=.9] (-11.4,1.05)--(-14.5,8.8) ;

\draw[ressort, decorate] (-19.7,2.5)--(-16,8.3) ;
\draw[ressort, decorate] (-11.4,1.05)--(-14.5,8.8) ;
\draw[ressort, decorate] (-19.5,1.8)--(-12.18,.2) ;

\draw[orange] [rotate around={160:(-14.7,3.6)}, line width=1.5pt, -stealth, dashed] (-10,2) arc (-30:210:4.5cm and 1cm);

\draw[dotted, thick] (-44,-3) -- (4,-3);

\draw[dashed, thick] (-15,-1.5) -- (-15,-5.5);
\draw[dashed, thick,->] (-15,-5.5) -- (-18,-7.5);
\draw[dashed, thick,->] (-15,-5.5) -- (-12,-7.5);

}}

\draw  (-26,12) node[diagonal fill={orange!40}{orange!40},
      text width=2cm, minimum height=1.1cm,
      text centered, rounded corners, draw, drop shadow] (0,0) {State-based};

\draw (-7,16.2) node{ \textbf{State} \textcolor{orange!80!black}{$\zeta$}};

    \path (5,5) pic[scale=.8]{molecule_full};

\draw  (-26,-2) node[diagonal fill={gray!40}{olive!40},
      text width=2.2cm, minimum height=1.1cm,
      text centered, rounded corners, draw, drop shadow] {Energy-based} ;

\node[circle, inner color = olive, outer color = white]
    at (-15,-8.2) {$\hspace{90pt}$};
    \path (-15,-12) pic{molecule_normal};

\draw (-15,-1) node{ \textbf{Energy} level \textcolor{olive!80!black}{$I$}};


\draw (0,-1) node{ \textbf{Energy} quantile $\mathcolor{gray!50!black}{q}$};

\node[circle,inner color = black, outer color = white] at (1.5,-11) {$\hspace{90pt}$};
\path (2.2,-14.2) pic{molecule_normal};

\path [fill=white, opacity=.25] (8,-3) rectangle (22,12);

\node[circle,inner color = black, outer color = white, opacity=.5]  at (1,-13.5+3.8)  {$\hspace{90pt}$};

\path (1.55,-13.65) pic{molecule_normal};

\path [fill=white, opacity=.25] (8,-3) rectangle (22,12);

\node[circle,inner color = black, outer color = white, opacity=.5] at (.5,-12.7+3.8) {$\hspace{90pt}$};

\path (.7,-12.8) pic{molecule_normal};

\path [fill=white, opacity=.25] (8,-3) rectangle (22,12);

\node[circle,inner color = black, outer color = white, opacity = .5] at (0,-12+3.8) {$\hspace{90pt}$};

\path (0,-12) pic{molecule_normal};

\end{tikzpicture}

    \caption{Schematic representation of the three main points of view for describing the internal structure of a polyatomic molecule.}
    \label{fig:modelling_choices}
\end{figure}
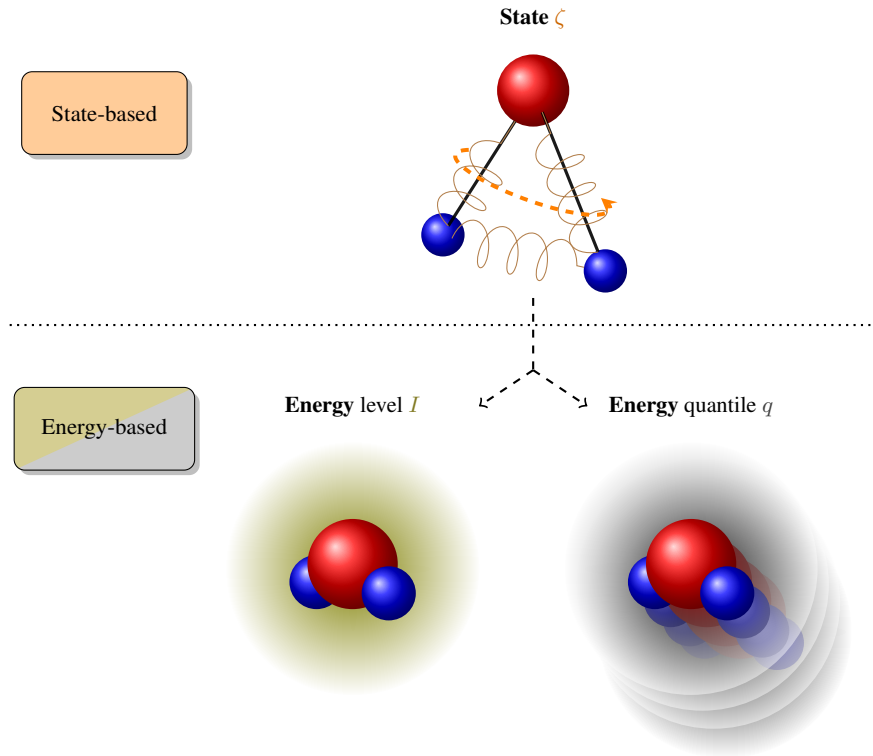

\noindent To sum up, there are three main points of view to describe the internal structure of a polyatomic molecule, a schematic representation of which is given in Fig.~\ref{fig:modelling_choices}.
A natural question to raise is then: \textbf{\emph{which one to choose between the three?}} As it turns out, all three have various characteristics that make them fit for different situations, as we detail in this sub-section.

\medskip

\noindent \emph{\textbf{State-based description.}} This is the \emph{simplest to build from physics} as any number of variables of various types can be considered. This is also the \emph{most complete} of all descriptions, as one can consider a model with all desired details, it therefore contains the ``full information'' (relatively to the underlying physical model) on the molecule. Finally, it is \emph{mathematically general} and, in fact, also includes the other descriptions. Nevertheless, \emph{the space of states could be large}, and \emph{neither the measure} $\mu$ \emph{nor the energy function} $\e$ \emph{are a priori trivial}\footnote{By ``trivial'', we mean, for $\mu$, to be a uniform measure, and for $\e$ to be the identity.}.

\smallskip

\noindent This description is particularly relevant in a \textbf{modelling context}, as well as for \textbf{general results}. 

\bigskip
\bigskip

\noindent \emph{\textbf{Energy-based descriptions.}} These are \emph{hardly constructed from physics, rather from a state-based description}. Here, \emph{the information is reduced to the energy}, and the underlying states are forgotten. The variable under consideration belongs to a \emph{one-dimensional space} (while in the state-based approach it could be multidimensional), and \emph{either the measure} $\mu$ \emph{or the energy function} $\e$ \emph{is trivial}.

\bigskip

\begin{minipage}{.04\textwidth}
    \phantom{x}
\end{minipage}
\begin{minipage}{.92\textwidth}
\emph{\textbf{Energy levels.}} The internal energy function is (almost) trivial: the energy associated with the energy level $I$ is $I + \e^0$. All information on the internal structure of the molecule is summarized in the measure $\mu_{\Bar{\e}}$. 

\smallskip

\noindent This model is particularly relevant in a context of \textbf{analysis}, typically when in-depth computations are involved. 
\end{minipage}

\bigskip

\begin{minipage}{.04\textwidth}
    \phantom{x}
\end{minipage}
\begin{minipage}{.92\textwidth}
\emph{\textbf{Energy quantiles.}} The measure on the space of quantiles is trivial: it is the Lebesgue measure. All information on the internal structure of the molecule is summarized in the quantile energy function $F^{\leftarrow}_{\mu_{\Bar{\e}}}$. 

\smallskip

\noindent This model is particularly relevant in a context of \textbf{particle-based numerical simulations}, both because of the \emph{uni-dimensionality} of the state space and the fact that it is endowed with a \emph{uniform measure}. This property is indeed fundamental in a particle-based simulation, as changing the state of a particle must not unintentionally also change its numerical weight.
\end{minipage}

\medskip

\noindent We sum up in Table~\ref{table:usefulwhen?} the characteristics and use of each description.

\begin{table}[h!] 
    \centering
  \def\arraystretch{1.8}
  \begin{tabularx}{\textwidth}{ 
  | >{\hsize=.6\hsize\linewidth=\hsize\centering}X 
  | >{\hsize=1.1\hsize\linewidth=\hsize\centering}X 
  || >{\hsize=1.1\hsize\linewidth=\hsize\centering}X 
  | >{\hsize=1.2\hsize\centering\arraybackslash}X |}
 \hline
\cellcolor{black!35} &  \cellcolor{white!90!orange} \textbf{\emph{State-based}: state} \textcolor{orange!80!black}{$\zeta$} & \cellcolor{white!90!olive} \textbf{\emph{Energy-based}: energy level} \textcolor{olive!80!black}{$I$}  &  \cellcolor{white!90!gray}  \textbf{\emph{Energy-based}: \hspace{-3pt}energy quantile} $\mathcolor{gray!70!black}{q}$  \\
        \hline 
  \cellcolor{black!5!white} \emph{Strengths} &  built from physics, general
  &     1-D space,
      trivial energy function
   &  1-D space, uniform measure  \\
       \hdashline 
     \cellcolor{black!5!white} \emph{Basic concepts} &  state space, measure and energy function &  measure & 
     energy function \\
        \hdashline 
\cellcolor{black!5!white}  \emph{Microscopic information} &   complete &  reduced to internal energy &  reduced to internal energy  \\
     \hline  
 \cellcolor{black!20!white!90!blue}\textbf{Adapted to}  & \cellcolor{blue!10!white} Modelling and general proofs &  \cellcolor{blue!10!white} Computations and technical proofs & \cellcolor{blue!10!white} Particle-based numerical simulations \\
\hline
\end{tabularx}
    \caption{Use of each description.}
    \label{table:usefulwhen?}
\end{table}

\newpage

\subsection{Model combination} \label{section:modelcomb}
\noindent Some molecular models can be written as \emph{independent combinations} of smaller ones. 

\subsubsection{Independent combination of two state-based models}

\noindent Consider two state-based models, one with a space of states $\E_1$ endowed with $\mu_1$ on which it is defined the energy function $\e_1$, and the other characterized by $\E_2$, $\mu_2$ and $\e_2$. We call \emph{independent combination of models $1$ and $2$} the internal state description defined by\footnote{Rigorously, $\mu_1$ and $\mu_2$ should be $\sigma$-finite.}

\smallskip

- a space of internal states $\E = \E_1 \times \E_2$ endowed with the measure $\mu_1 \otimes \mu_2$,

- a \emph{total} internal energy function $\e : \begin{cases}
    \E_1 \times \E_2 \quad  \to \quad \R \\
    (\zeta_1,\zeta_2) \; \;\mapsto \; \; \e_1(\zeta_1) + \e_2 (\zeta_2)
\end{cases}$.

\noindent We mention right away that the \emph{number of internal degrees of freedom} (later defined in Section~\ref{section3}) associated with the combined model is the \textbf{sum} of the numbers of internal degrees of freedom associated with model 1 and model 2, that is $\delta_{1 \& 2} = \delta_1 + \delta_2$.

\subsubsection{\emph{Separated}- and \emph{total}-energy-based corresponding descriptions}

\paragraph{\phantom{abc} \textbf{\emph{Separated}-internal energy levels}} As in the combination model described above there are two types of energies, we may construct a \emph{separated-internal energy levels} model from the independent combination of models $1$ and $2$ as follows:

\smallskip

\noindent The internal variable is the \emph{separated-internal energy levels} $(I_1,I_2) \in \R_+^2$, its associated energy is $I_1 + I_2 + \e^0_1 + \e^0_2$. The separated-internal energy levels space $\R_+^2$ is endowed with the measure $\mu_{\bar{\e}_1} \otimes \mu_{\bar{\e}_2}$, where $\mu_{\bar{\e}_i}$ is the energy law of the $i^{th}$ model, $\mu_{\bar{\e}_i} = \mu_i \# \bar{\e}_i$. We also recall that $I_i$ relates to $\bar{\e}_i(\zeta_i)$ (see Subsubsection~\ref{section:energy_levels}). 

\medskip

\paragraph{\phantom{abc} \textbf{\emph{Separated}-internal energy quantiles}} Similarly, we obtain a \emph{separated-internal energy quantiles} as follows:
the internal variable is the \emph{separated-internal energy quantiles} $(q_1,q_2) \in (0,q^{\max}_1) \times (0,q^{\max}_2)$, its associated energy is $F^{\leftarrow}_{\mu_{\bar{\e}_1}}(q_1) + F^{\leftarrow}_{\mu_{\bar{\e}_2}}(q_2) + \e^0_1 + \e^0_2$, where $F^{\leftarrow}_{\mu_{\bar{\e}_i}}$ is the internal energy quantile function associated with the $i^{th}$ model (see Subsubsection~\ref{section:energy_quantiles}). The separated-internal energy levels space $(0,q^{\max}_1) \times (0,q^{\max}_2)$ is endowed with the Lebesgue measure. 

\paragraph{\phantom{abc} \textbf{Link between \emph{separated}- and \emph{total}-internal energy levels descriptions}}

Let us recall that the internal energy levels description associated with the (total) combined model is given as follows:

\smallskip

\noindent The internal variable is the total internal energy level $I \in \R_+$, the latter space being endowed with the measure $\mu_{\bar{\e}} = \mu \# \bar{\e}$, and the energy associated with the level $I$ is $I + \e^0$.

\medskip

\noindent We can relate the \emph{separated}- and \emph{total}-internal energy levels descriptions through the following relationships:

\begin{align}
    I &= I_1 + I_2, \\
    \e^0 &= \e^0_1 + \e^0_2, \\
    \mu_{\bar{\e}} &= \mu_{\bar{\e}_1} \ast \mu_{\bar{\e}_2}, \label{eq:conv}
\end{align}
where $\ast$ stands for the convolution of measures.

\medskip

\noindent Equation~\eqref{eq:conv} should be understood through the probability setting interpretation of Subsection~\ref{section:proba}: the total internal energy function $\e$, seen as a ``random variable'', is the sum of the two \emph{independent} ``random variables'' $\e_1$ and $\e_2$. As such, the ``law'' of $\e$ is indeed given by the convolution of the ``laws'' of $\e_1$ and $\e_2$.

\section{A typical model for a diatomic molecule} \label{section2}
In order to provide an example typically considered in the kinetic theory literature,
let us focus, as in the diagram below, on a \emph{diatomic molecule} (composed of two atoms), for which we want to take into account both \emph{rotation} and \emph{vibration} phenomena.

\subsection{Molecular model and internal states description}

In this example, we choose 
to describe the rotation via a \emph{classical} approach, by considering the \emph{angular velocity}\footnote{The angular velocity has only two components here because the molecule is linear.} $\omega \in \R^2$, and the vibration via a quantum harmonic oscillator model,

\medskip

\begin{minipage}{.27\textwidth}
   \hspace{-1em} \begin{tikzpicture}[scale=.3,decoration={coil,aspect=0.4,segment length = 3mm, amplitude=3mm}]
\tikzset{ressort/.style={thick,gray,smooth}}
\shade[ball color=green!50!red] (-5,2) circle (1.2);
\shade[ball color=green!50!red] (4,0) circle (1.2);
\draw[ultra thick] (-4.5,1.8)--(2.82,.2) ;

\draw (-.5,3) node{\AxisRotator[rotate=-100]};
\draw[->,thick] (-.5,3)--(-.1,5);
\draw (-1,5) node{$\omega$};
\draw[ressort,decorate] (-4.5,1.8)--(2.82,.2) ;
\draw (-.4,-1.2) node{$n$};
\end{tikzpicture}
\end{minipage}
\begin{minipage}{.68\textwidth}
       denoting by $n \in \N$ the vibration mode. With a rigid-rotor assumption, making rotation and vibration independent, the internal energy associated with the $(\omega,n)$ state can then be written as the sum of the kinetic energy of rotation and the energy of vibration:
      \[
\e(\omega,n) = \frac12 \J |\omega|^2 + \left( n  + \frac12 \right) \Delta \epsilon,
\]
where \(\J > 0\) is the moment of inertia and $\Delta \epsilon > 0$ the energy gap between two vibration levels.
\end{minipage}

\bigskip

\noindent Let us provide the internal (state-based, energy-based) descriptions that correspond to the above model.

\subsubsection*{Corresponding \emph{internal states} description.} The internal states description is trivial to construct from the molecular model:

\medskip

-- the space of internal states is $\E = \R^2 \times \N$ endowed with the measure $\mu = \mathrm{Lebesgue}_{\R^2} \otimes \mathrm{Counting}_{\N}$ (as there is no degeneracy)

-- the internal energy function is
\( \displaystyle
\e : \begin{cases}
    \R^2 \times \N \; &\longrightarrow \qquad \R \\
    (\omega, n) \; &\mapsto \; \frac12 \J |\omega|^2 + \left( n  + \frac12 \right) \Delta \epsilon,
\end{cases}
\)
which minimum is \(\e^0 = \frac12 \Delta \epsilon\).

\smallskip

\noindent Note that this model has the shape of an \emph{independent combination of two sub-models} (one for rotation, one for vibration), see Subsection~\ref{section:modelcomb}. 

\subsection{Corresponding \emph{total} internal energy-based descriptions}

Let us provide the corresponding \emph{total} internal energy-based descriptions (see Subsection~\ref{section:modelcomb}).

\subsubsection{\emph{Total internal energy levels} description.} We obtain the total internal energy levels description by computing the \emph{energy law} $\mu_{\bar{\varepsilon}}$. It is provided, for any $I \geq 0$, by
\[
\mu_{\bar{\varepsilon}}([0,I)) = \mu (\{\e - \e^0 < I\}) = \sum_{n \in \N} \int_{\R^2} \mathds{1}_{\{\frac12 \J |\omega|^2 + n \Delta \epsilon < I\}} \, \dd \omega,
\]
yielding
\begin{equation}
    \mu_{\bar{\varepsilon}}([0,I)) = \int_0^I \varphi_{\mathrm{total}}(I_*) \, \dd I_*,
\end{equation}
with
\begin{equation} \label{eq:varphi_diatomic}
   \varphi_{\mathrm{total}}(I) = \frac{2\pi}{\J} \left\lceil \frac{I}{\Delta \epsilon} \right\rceil,
\end{equation}
where \(\left\lceil \cdot \right\rceil\) stands for the upper integer part. The internal energy levels description is therefore composed of:

\medskip

-- the space of total internal energy levels $\R_+$ endowed with the measure $\varphi_{\mathrm{total}}(I) \, \dd I$, with $\varphi_{\mathrm{total}}$ given by~\eqref{eq:varphi_diatomic}. It is plotted in Fig.~\ref{figure:energylaws}.

-- the fundamental energy level $\e^0 = \frac12 \Delta \epsilon$.

\subsubsection{\emph{Total internal energy quantiles} description.} \noindent We obtain the total internal energy quantiles description by computing the \emph{energy quantile function} $F^{\leftarrow}_{\mu_{\bar{\varepsilon}}}$. First note that
\[q^{\max} = \mu(\E) = + \infty.\]
The quantile function is such that for any $q \in (0, +\infty)$ and $I \geq 0$,
\[
F^{\leftarrow}_{\mu_{\bar{\varepsilon}}}(q) = \inf \left\{ I \geq 0, \; q \leq \mu_{\bar{\varepsilon}}([0,I)) \right\} \quad = \quad \inf \left\{ I \geq 0, \; \; \; q \leq \int_{0}^I \frac{2\pi}{\J} \left\lceil \frac{I_*}{\Delta \epsilon} \right\rceil \, \dd I_* \right\},
\]
from which one finds
\begin{equation} \label{eq:quantilefunction_diatomic}
 \forall \, q \geq 0, \qquad   F^{\leftarrow}_{\mu_{\bar{\varepsilon}}}(q) = \Delta \epsilon \left(\frac{\hat{q}}{\ell(\hat{q})+1} + \frac{\ell(\hat{q})}{2} \right),
\end{equation}
with
\[
\hat{q} := \frac{\J}{2\pi \; \Delta \epsilon} \; q, \qquad \ell(\hat{q}) := \left\lfloor \frac{4 \, \hat{q}}{\sqrt{1 + 8 \, \hat{q}} \, + 1} \right\rfloor.
\]
The internal energy quantiles description is therefore composed of:

\medskip

-- the space of total internal energy quantiles $(0,q^{\max}) = (0, + \infty)$.

-- the total internal energy quantile function $F^{\leftarrow}_{\mu_{\bar{\varepsilon}}}$ given by~\eqref{eq:quantilefunction_diatomic} and the fundamental energy level $\e^0 = \frac12 \Delta \epsilon$.

\subsection{Corresponding \emph{separated} internal energy-based descriptions}

We now provide the corresponding \emph{separated} internal energy-based descriptions. We highlight again for that matter that the internal states model we consider here is an \emph{independent combination} of a model for rotation and a model for vibration:

\medskip

-- the space of internal states is $\E = \E_{rot} \times \E_{vib}$  endowed with the measure $\mu = \mu_{rot} \otimes \mu_{vib}$, where $\E_{rot} = \R^2$, $\E_{vib} = \N$, $\mu_{rot} = \mathrm{Lebesgue}_{\R^2}$ and $\mu_{vib} = \mathrm{Counting}_{\N}$.

\medskip

-- the internal energy function is \(\e(\omega,n) = \e_{rot}(\omega) + \e_{vib}(n)\), with \\
\[\e_{rot} : \begin{cases}
    \R^2\; &\longrightarrow \qquad \R \\
    \omega \; &\mapsto \; \frac12 \J |\omega|^2
\end{cases} \quad \text{and} \quad  \e_{vib} : \begin{cases}
    \N \; &\longrightarrow \qquad \R \\
    n \; &\mapsto \; \left( n  + \frac12 \right) \Delta \epsilon
\end{cases}.
\]

\smallskip

\noindent We have in particular $\e_{rot}^0 = 0$ and $\e_{vib}^0 = \frac12 \Delta \epsilon$.

\medskip

\noindent We may then consider two energy variables, one associated with rotation and one with vibration.

\subsubsection{\emph{Separated ro-vib internal energy levels} description.} We obtain the \emph{separated ro-vib} internal energy levels description by computing the \emph{energy laws} $\mu_{\bar{\varepsilon}_{rot}}$ and $\mu_{\bar{\varepsilon}_{vib}}$, which are, for any $I > 0$, such that 
\begin{align*}
\mu_{\bar{\varepsilon}_{rot}}([0,I)) &= \mu_{rot} (\{\e_{rot} - \e_{rot}^0 < I\}) = \int_{\R^2} \mathds{1}_{\{\frac12 \J |\omega|^2 < I\}} \, \dd \omega, \\
 \mu_{\bar{\varepsilon}_{vib}}([0,I)) &= \mu_{vib} (\{\e_{vib} - \e_{vib}^0 < I\}) = \sum_{n \in \N}\mathds{1}_{\{ n \Delta \epsilon < I\}},
\end{align*}

yielding
\begin{equation} \label{eq:energylaws_rovib}
    \dd \mu_{\bar{\varepsilon}_{rot}}(I) = \varphi_{rot}(I) \, \dd I, \qquad \mu_{\bar{\varepsilon}_{vib}} = \sum_{n \in \N} \boldsymbol{\delta}_{n \Delta \epsilon},
\end{equation}
where $\varphi_{rot}(I) := \frac{2 \pi}{\mathcal{J}}$ is therefore constant, and $\boldsymbol{\delta}_{n \Delta \epsilon}$ stands for the Dirac mass at $n \Delta \epsilon$. As a remark, notice that the \emph{total} internal energy law $\mu_{\bar{\varepsilon}}$ indeed satisfies
\(\mu_{\bar{\varepsilon}} = \mu_{\bar{\varepsilon}_{rot}} \ast \mu_{\bar{\varepsilon}_{vib}}\). The \emph{separated ro-vib} internal energy levels description is therefore composed of:

\medskip

-- the space of \emph{separated ro-vib} internal energy levels $\R_+^2$ endowed with the measure $\mu_{\bar{\varepsilon}_{rot}} \otimes \mu_{\bar{\varepsilon}_{vib}}$, each given in~\eqref{eq:energylaws_rovib}.

\smallskip

-- the fundamental energy level $\e_{rot}^0 + \e_{vib}^0= \frac12 \Delta \epsilon$.

\medskip

\noindent As it is equivalent to consider $\R_+$ endowed with a measure that is a countable sum of Dirac masses, and $\N$ endowed with a counting measure, we may prefer here to write the \emph{separated ro-vib} internal energy levels description as:

\medskip

-- the space of \emph{separated ro-vib} internal energy levels is $\R_+ \times \N$ endowed with the measure $\mu_{\bar{\varepsilon}_{rot}} \otimes \mathrm{Counting}_{\N}$, with $\mu_{\bar{\varepsilon}_{rot}}$ given in~\eqref{eq:energylaws_rovib}.

\smallskip

-- the internal energy associated with $(I_{rot},n_{vib}) \in \R_+ \times \N$ is $I_{rot} + \left(n + \frac12 \right) \Delta \epsilon$.

\subsubsection{\emph{Separated ro-vib internal energy quantiles} description.} \noindent We obtain the \emph{separated ro-vib} internal energy quantiles description by computing both rotational and vibrational energy quantile functions $F^{\leftarrow}_{\mu_{\bar{\varepsilon}_{rot}}}$ and $F^{\leftarrow}_{\mu_{\bar{\varepsilon}_{vib}}}$. We have
\[q^{\max}_{rot} = \mu_{rot}(\E_{rot}) = + \infty, \qquad q^{\max}_{vib} = \mu_{vib}(\E_{vib}) = + \infty.\]
The rotational quantile function is such that for any $q \in (0, +\infty)$ and $I \geq 0$,
\[
F^{\leftarrow}_{\mu_{\bar{\varepsilon}_{rot}}}(q) = \inf \left\{ I \geq 0, \; q \leq \mu_{\bar{\varepsilon}_{rot}}([0,I)) \right\} \quad = \quad \inf \left\{ I \geq 0, \; \; \; q \leq \frac{2\pi}{\J} I \right\},
\]
from which one finds
\begin{equation} \label{eq:quantilefunction_rot}
 \forall \, q > 0, \qquad   F^{\leftarrow}_{\mu_{\bar{\varepsilon}_{rot}}}(q) = \frac{\J}{2\pi} q.
\end{equation}
Similarly,
\begin{equation} \label{eq:quantilefunction_vib}
 \forall \, q > 0, \qquad   F^{\leftarrow}_{\mu_{\bar{\varepsilon}_{vib}}}(q) = \Delta \epsilon \left\lfloor q \right\rfloor.
\end{equation}

The \emph{ro-vib} internal energy quantiles description is therefore composed of:

\medskip

-- the space of \emph{ro-vib} internal energy quantiles $(0, + \infty)^2$.

-- the \emph{ro-vib} internal energy quantile function $(q_{rot},q_{vib}) \mapsto F^{\leftarrow}_{\mu_{\bar{\varepsilon}_{rot}}}(q_{rot}) + F^{\leftarrow}_{\mu_{\bar{\varepsilon}_{vib}}}(q_{vib})$, which functions are given in~\eqref{eq:quantilefunction_rot}--\eqref{eq:quantilefunction_vib}, and the fundamental energy level $\e^0 = \frac12 \Delta \epsilon$.

\subsection{Plots of the energy laws and quantile functions}

\begin{figure}[h!] 
\centering
\begin{subfigure}[t!]{0.48\textwidth}    \includegraphics[width=\linewidth]{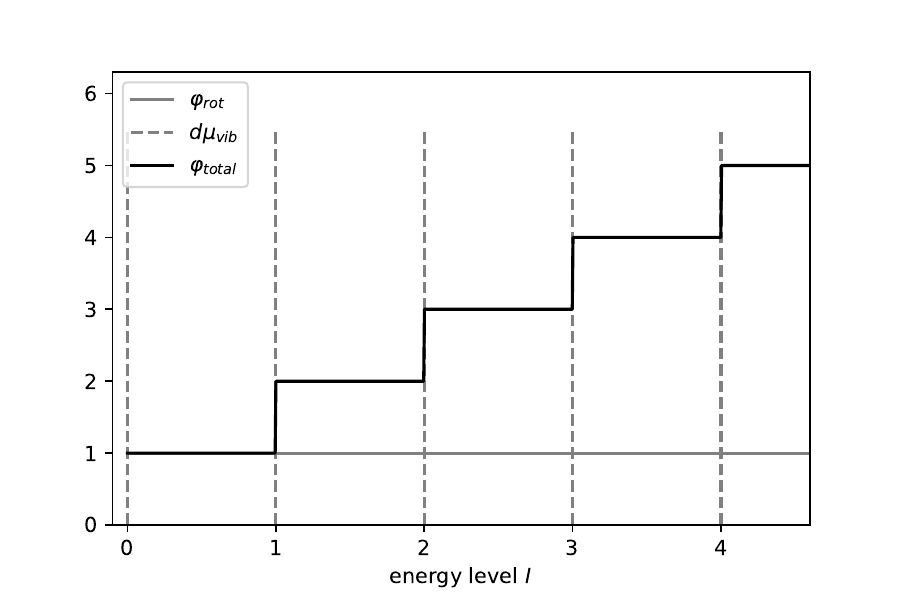}
    \caption{Plot of the energy law densities for the cases of rotation, vibration, and total, respectively, with $x$-axis rescaled by $\Delta \epsilon$ and $y$-axis rescaled by $\frac{\mathcal{J}}{2 \pi \, \Delta \epsilon}$.}
    \label{figure:energylaws}
    \end{subfigure} \hspace{5pt}
\centering
\begin{subfigure}[t!]{0.48\linewidth}
\includegraphics[width=\linewidth]{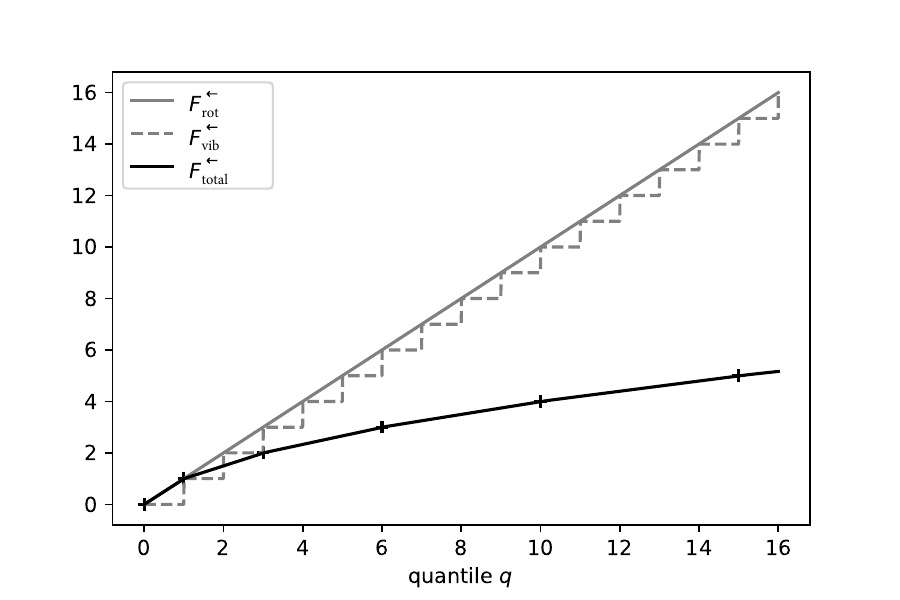}
    \caption{Plot of the quantile functions for the cases of rotation, vibration, and total,  respectively, with $x$-axis rescaled by $\frac{2 \pi \, \Delta \epsilon}{\mathcal{J}}$ and $y$-axis rescaled by $\frac{1}{\Delta \epsilon}$.}
    \label{figure:quantilefunctions}
\end{subfigure}
    \caption{Plots of energy laws and quantile functions associated with the discussed model.}
\end{figure}

\noindent As a matter of illustration, we plot in Fig.~\ref{figure:energylaws}--\ref{figure:quantilefunctions} the internal energy laws and quantile functions computed earlier in this sub-section, given in Equations~\eqref{eq:varphi_diatomic}--\eqref{eq:quantilefunction_vib}. Vertical lines in Fig.~\ref{figure:energylaws} represent Dirac masses. In Fig.~\ref{figure:quantilefunctions}, $F^{\leftarrow}_{\text{total}}$ (which is defined by~\eqref{eq:quantilefunction_diatomic}) is piecewise-affine, with breaks happening at every $q$ such that there is an integer $n$ such that $q = \frac{n (n+1)}{2}$.

\section{Resulting models at mesoscopic level} \label{section3}

In the previous section, we focused on mathematical settings for the internal description of polyatomic molecules. In this section, we show how this translates at the mesoscopic level, as well as in which situation one is allowed, or not, to navigate between the state-based and energy-based descriptions.

\smallskip

\noindent Consider a gas composed of polyatomic molecules of the same species, which internal state model is given by a space of internal states $\E$ endowed with a measure $\mu$ and an internal energy function $\e$. We denote by $\mu_{\bar{\e}}$ its associated internal energy law and $F^{\leftarrow}_{\mu_{\bar{\e}}}$ its internal energy quantile function (see Subsubsections~\ref{section:energy_levels} and~\ref{section:energy_quantiles}). 

\medskip

\noindent Each molecule can then be described at each time by its position, velocity, and internal variable (state, energy level, or energy quantile). The mathematical object of study is the molecular density function $f$, depending on time $t$, position $x$, velocity $v$ and internal variable.

\subsection{Macroscopic quantities}

\noindent As a matter of illustration, we provide in Table~\ref{table:macro} the formulas to obtain the density, bulk velocity and energy for each internal description. For clarity, the macroscopic variables $t$ and $x$ are left aside. The temperature $T$ is related to the energy $E$ via a potentially highly non-linear relationship (see~\cite[Eqs. (3.2) and (4.8)]{GF}).

\begin{table}[htbp] 
    \centering
  \def\arraystretch{2}
  \begin{tabularx}{\textwidth}{
  | >{\hsize=.04\hsize\linewidth=\hsize\centering}X
  | >{\hsize=1.27\hsize\linewidth=\hsize\centering}X 
  | >{\hsize=1.27\hsize\linewidth=\hsize\centering}X 
  | >{\hsize=1.42\hsize\linewidth=\hsize\centering\arraybackslash}X |}
 \hline
\rowcolor{black!15} \cellcolor{black!35} &  \cellcolor{orange!15} \textbf{States}& \cellcolor{olive!15} \textbf{Energy \emph{levels}} &  \cellcolor{black!10} \textbf{Energy \emph{quantiles}} \\
\hline
    \cellcolor{black!15}  \vspace{1em}  $f$  \vspace{1em}  & $f(v,\textcolor{orange!90!black}{\zeta})$ & $f(v,\textcolor{olive}{I})$ & $f(v,\textcolor{black!60!white}{q})$\\
\hdashline
\cellcolor{black!15} \vspace{1em} $\rho$ \vspace{1em} &   $\displaystyle \hspace{-1em} \iint_{\R^3 \times \textcolor{orange!90!black}{\E}} f(v,\textcolor{orange!90!black}{\zeta}) \, \dd v \, \dd \textcolor{orange!90!black}{\mu}(\textcolor{orange!90!black}{\zeta})$ & $\displaystyle \hspace{-1em} \iint_{\R^3 \times \textcolor{olive}{\R_+}} \hspace{-13pt} f(v,\textcolor{olive}{I}) \, \dd v \, \dd \textcolor{olive}{\mu_{\bar{\e}}}(\textcolor{olive}{I})$ & $\displaystyle \hspace{-1em} \iint_{\R^3 \times \textcolor{black!60!white}{(0,q^{\max})}} \hspace{-6pt} f(v,\textcolor{black!60!white}{q}) \, \dd v \, \dd \textcolor{black!60!white}{q}$ \\
\hdashline
\cellcolor{black!15} \vspace{1em} $u$ \vspace{1em} & $\displaystyle \hspace{-7pt} \hspace{-.5em} \frac{1}{\rho} \iint_{\R^3 \times \textcolor{orange!90!black}{\E}}\hspace{-1em} v \, f(v,\textcolor{orange!90!black}{\zeta}) \, \dd v \, \dd \textcolor{orange!90!black}{\mu}(\textcolor{orange!90!black}{\zeta})$ & $\displaystyle \hspace{-5pt} \hspace{-1em} \frac{1}{\rho} \iint_{\R^3 \times \textcolor{olive}{\R_+}} \hspace{-15pt} vf(v,\textcolor{olive}{I}) \, \dd v \, \dd \textcolor{olive}{\mu_{\bar{\e}}}(\textcolor{olive}{I})$ & $\displaystyle \hspace{-1em} \frac{1}{\rho} \iint_{\R^3 \times \textcolor{black!60!white}{(0,q^{\max})}} \hspace{-6pt} \hspace{-1em} vf(v,\textcolor{black!60!white}{q}) \, \dd v \, \dd \textcolor{black!60!white}{q}$ \\
\hdashline
\cellcolor{black!15} \vspace{1em} $E$ \vspace{1em} &  
$\underset{ \displaystyle \times \, f(v,\textcolor{orange!90!black}{\zeta}) \, \dd v \, \dd \textcolor{orange!90!black}{\mu}(\textcolor{orange!90!black}{\zeta})}{\displaystyle \hspace{-7pt}  \underset{\R^3 \times \textcolor{orange!90!black}{\E}}{\iint \hspace{20pt}}  \hspace{-20pt}\left(\frac12 |v-u|^2 + \textcolor{orange!90!black}{\e}(\textcolor{orange!90!black}{\zeta}) \right) }$ & $\underset{\displaystyle  \times \, f(v,\textcolor{olive}{I}) \, \dd v \, \dd \textcolor{olive}{\mu_{\bar{\e}}}(\textcolor{olive}{I})}{\displaystyle \hspace{-12pt} \underset{\R^3 \times \textcolor{olive}{\R_+}}{\iint \hspace{20pt}} \hspace{-20pt}\left(\frac12 |v-u|^2 + \textcolor{olive}{I} + \e^0 \right) }$ & $\underset{\displaystyle  \times \, f(v,\textcolor{black!60!white}{q}) \, \dd v \, \dd \textcolor{black!60!white}{q}}{\displaystyle \hspace{-16pt} \underset{\R^3 \times \textcolor{black!60!white}{(0,q^{\max})}}{\iint \hspace{35pt}} \hspace{-40pt}  \left(\frac12 |v-u|^2 + \textcolor{black!60!white}{F^{\leftarrow}_{\mu_{\bar{\e}}}}(\textcolor{black!60!white}{q}) + \e^0 \right) }$ \\
\hline
\end{tabularx}
    \caption{Standard macroscopic quantities for each internal description, with $f$ the mesoscopic molecular density, $\rho$ the macroscopic molecular density, $u$ the bulk velocity and $E$ the energy.}
    \label{table:macro}
\end{table}

\subsection{Equivalence of the state-based and energy-based approaches at the mesoscopic level}

The energy-based internal descriptions contain less information on the molecule than the full state-based one. Does this \emph{loss of information} at the microscopic level (the fact of forgetting the internal states in order to keep only the distribution of energy levels) have an impact on the model at the \emph{mesoscopic} level? 

\medskip

\noindent The answer is that \emph{at thermodynamical equilibrium, the approaches are always equivalent}, and \emph{out of thermodynamical equilibrium, it depends on the molecular interactions}.

\medskip

\noindent Typically, in a Boltzmann context, the state-based and energy-based Boltzmann models would be equivalent if, and only if, the \emph{cross-section}, which characterizes the distribution of post-collisional states as a function of pre-collisional ones, depends on the internal states only via their energy.

\medskip

\noindent Let us take as an example the model describing the rotation and vibration of a diatomic molecule presented in Section~\ref{section2}, and consider a Boltzmann model associated with this molecular model. If the cross-section depends, as far as the internal states are concerned, only \emph{separately} on the energy of rotation and the energy of vibration, then the internal information can be reduced to the \emph{separate} internal energies, but not to the \emph{total} internal energy, without loss of information at kinetic level. Then, the state-based and  \emph{separated}-energy-based Boltzmann models are equivalent, but \emph{not} the \emph{total}-energy-based one. On the other hand, if the cross-section depends, as far as the internal states are concerned, only on the \emph{total energy}, that is the sum of the energy of rotation and the energy of vibration, then the internal information can be reduced to the \emph{total} internal energy without loss of information at kinetic level. Then, the state-based, \emph{separated}-energy-based and \emph{total}-energy-based Boltzmann models are all equivalent.

\subsection{Equilibrium quantities} \label{subsec:equilqties}

We conclude this note with the mention of the form of the equilibrium distribution in a \emph{classical polyatomic context}, the Maxwellian distribution, for each approach, as well as two major equilibrium quantities that are the \emph{number of degrees of freedom} and the \emph{heat capacity at constant volume}.

\subsubsection{Maxwellian measure}
The Maxwellian measure is related to collision equilibrium, and we provide its form for each approach in Table~\ref{table:maxwellian}. The function $M \equiv M[\rho,u,T]$ stands for the monoatomic Maxwellian distribution defined, for any $v \in \R^3$, by
\begin{equation}
    M(v) := \rho\, (2 \pi k_B T)^{-\frac32} \, e^{-\frac{|v-u|^2}{2 \, k_B T}},
\end{equation}
where $k_B$ stands for the Boltzmann constant, while $Z$ denotes the internal partition function, defined, for any $\beta>0$, by,
\begin{equation}
    Z(\beta) := \int_{\textcolor{orange!90!black}{\E}} e^{-\beta \, \textcolor{orange!90!black}{\bar{\e}}(\textcolor{orange!90!black}{\zeta})} \, \mathrm{d}\textcolor{orange!90!black}{\mu}(\textcolor{orange!90!black}{\zeta}) = \int_{\textcolor{olive}{\R_+}} e^{-\beta \, \textcolor{olive}{I}} \, \mathrm{d}  \textcolor{olive}{\mu_{\bar{\e}}}(\textcolor{olive}{I}) = \int_{\textcolor{black!60!white}{0}}^{\textcolor{black!60!white}{q^{\max}}} e^{-\beta \, \textcolor{black!60!white}{F^{\leftarrow}_{\mu_{\bar{\e}}}}(\textcolor{black!60!white}{q})} \, \mathrm{d} \textcolor{black!60!white}{q}.
\end{equation}

\newpage

\begin{table}[htbp]
    \centering
    \def\arraystretch{2.2}
    \begin{tabularx}{\textwidth}{|>{\centering\arraybackslash}X|>{\centering\arraybackslash}X|>{\centering\arraybackslash}X|}
        \hline
        \cellcolor{orange!15} \textbf{States}, on $\R^3 \times \textcolor{orange!90!black}{\E}$ & \cellcolor{olive!15} \textbf{Energy levels}, on $\R^3 \times \textcolor{olive}{\R_+}$ & \cellcolor{black!10} \textbf{Energy quantiles}, on $\R^3 \times \textcolor{black!60!white}{(0,q^{\max})}$ \\
        \hline
        $\underset{\displaystyle \phantom{\int} \times e^{-\frac{\textcolor{orange!90!black}{\bar{\e}}(\textcolor{orange!90!black}{\zeta})}{k_B T}} \,\mathrm{d}v \, \mathrm{d} \textcolor{orange!90!black}{\mu}(\textcolor{orange!90!black}{\zeta})}{\displaystyle
       M(v) \, Z\left(\frac{1}{k_B T}\right)^{-1}_{\phantom{\displaystyle \int}}}$ & 
        $\underset{\displaystyle \phantom{\int} \times e^{-\frac{\textcolor{olive}{I}}{k_B T}} \,\mathrm{d}v \, \mathrm{d} \textcolor{olive}{\mu_{\bar{\e}}}(\textcolor{olive}{I})}{\displaystyle
       M(v) \, Z\left(\frac{1}{k_B T}\right)^{-1}_{\phantom{\displaystyle \int}}}$ & 
        $\underset{\displaystyle \phantom{\int} \times e^{-\frac{\textcolor{black!60!white}{F^{\leftarrow}_{\mu_{\bar{\e}}}}(\textcolor{black!60!white}{q})}{k_B T}} \, \mathrm{d}v \, \mathrm{d} \textcolor{black!60!white}{q}}{\displaystyle
       M(v) \, Z\left(\frac{1}{k_B T}\right)^{-1}_{\phantom{\displaystyle \int}}}$ \\
        \hline
    \end{tabularx}
    \caption{Maxwellian measure on $\R^3 \hspace{-1pt} \times \cdot  $ for each approach, for some fixed $\rho, u$ and $T$.}
    \label{table:maxwellian}
\end{table}
\vspace{-.5em}
\noindent Finally, we provide the formulas for the \textbf{total number of degrees of freedom} $3 + \delta(T)$ ($3$ is the number of \emph{translational} degrees of freedom and $\delta$ is called number of \emph{internal} degrees of freedom), depending on the temperature of the system $T>0$, by
\vspace{-1.3em}
\begin{equation} \label{eq:delta1}
\begin{split}
3 + \delta(T) := 3 + 2 \beta \,  \frac{\displaystyle \int_{\textcolor{orange!90!black}{\E}} \textcolor{orange!90!black}{\bar{\e}}(\textcolor{orange!90!black}{\zeta}) \, e^{-\beta \, \textcolor{orange!90!black}{\bar{\e}}(\textcolor{orange!90!black}{\zeta})} \, \mathrm{d} \, \textcolor{orange!90!black}{\mu}(\textcolor{orange!90!black}{\zeta}) }{\int_{\textcolor{orange!90!black}{\E}} e^{-\beta \, \textcolor{orange!90!black}{\bar{\e}}(\textcolor{orange!90!black}{\zeta})} \, \mathrm{d} \, \textcolor{orange!90!black}{\mu}(\textcolor{orange!90!black}{\zeta}) }  =  3 +  2 \beta \,  \frac{\displaystyle  \int_{\textcolor{olive}{\R_+}} \textcolor{olive}{I} \, e^{-\beta \, \textcolor{olive}{I}} \, \mathrm{d}  \textcolor{olive}{\mu_{\bar{\e}}}(\textcolor{olive}{I})}{\displaystyle  \int_{\textcolor{olive}{\R_+}} e^{-\beta \, \textcolor{olive}{I}} \, \mathrm{d}  \textcolor{olive}{\mu_{\bar{\e}}}(\textcolor{olive}{I})} \\
= 3 +  2 \beta \,  \frac{\displaystyle \int_{\textcolor{black!60!white}{(0,q^{\max})}} \textcolor{black!60!white}{F^{\leftarrow}_{\mu_{\bar{\e}}}}(\textcolor{black!60!white}{q}) \, e^{-\beta \, \textcolor{black!60!white}{F^{\leftarrow}_{\mu_{\bar{\e}}}}(\textcolor{black!60!white}{q})} \, \mathrm{d} \textcolor{black!60!white}{q}}{\int_{\textcolor{black!60!white}{(0,q^{\max})}} e^{-\beta \, \textcolor{black!60!white}{F^{\leftarrow}_{\mu_{\bar{\e}}}}(\textcolor{black!60!white}{q})} \, \mathrm{d} \textcolor{black!60!white}{q}},
\end{split}
\end{equation}
where we denoted $\beta = \frac{1}{k_B T}$, or equivalently
\begin{equation} \label{eq:delta2}
3 + \delta(T) := 3 -  2 \beta \,  (\log Z)' \left( \beta \right).
\end{equation}
The total number of degrees of freedom can also be interpreted as the \emph{expectation} of the random variable $(v,\zeta) \mapsto \frac{2}{k_B T} \left(\frac12 |v-u|^2 + \bar{\e}(\zeta) \right)$ on the space $\R^3 \times \E$ endowed with the Maxwellian probability (the measure given in Table~\ref{table:maxwellian} with $\rho=1$).

\medskip

\noindent As for the \textbf{heat capacity at constant volume} $c_V(T)$, also depending on the temperature $T > 0$ it is given by
\begin{equation}
    c_V(T) := \frac32 + \frac12 \, \frac{\mathrm{d} (T \delta(T))}{\mathrm{d} T}.
\end{equation}
Other formulas for $c_V$ can be straightforwardly deduced from~\eqref{eq:delta1}--\eqref{eq:delta2}. It can also be interpreted as the \emph{variance} of the random variable $(v,\zeta) \mapsto \frac{1}{k_B T} \left(\frac12 |v-u|^2 + \bar{\e}(\zeta) \right)$ on the space $\R^3 \times \E$ endowed with the Maxwellian probability.

\medskip


As an example, for the diatomic molecular model presented in Sec. \ref{section2}, the total number of degrees of freedom writes, as a function of the temperature $T > 0$,
\begin{equation} \label{eq:delta}
    3 + \delta (T) = \underbracket{\, 3 \, }_{\text{translation}} + \underbracket{\, 2 \, }_{\text{rotation}} +  \underbracket{\, \displaystyle 2 \left( \frac{\Delta \epsilon}{k_B T} \right) \; \frac{1}{\displaystyle e^{\frac{\Delta \epsilon}{k_B T}}-1} \, }_{\text{vibration}},
\end{equation}
while the heat capacity at constant volume writes
\begin{equation} \label{eq:cv}
    c_V(T) = \underbracket{\, \frac{3}{2}  \, }_{\text{translation}} + \underbracket{\, 1 \, }_{\text{rotation}} +  \underbracket{\, \left(\frac{\Delta \epsilon}{k_B T} \right)^2 \frac{e^{\frac{\Delta \epsilon}{k_B T}}}{(e^{\frac{\Delta \epsilon}{k_B T}}-1)^2} \, }_{\text{vibration}},
\end{equation}
where we recall that $\Delta \epsilon$ is the energy gap between two vibration modes.

\medskip

\begin{figure}[h!] 
\centering
\begin{subfigure}[t!]{0.48\textwidth}    \includegraphics[width=\linewidth]{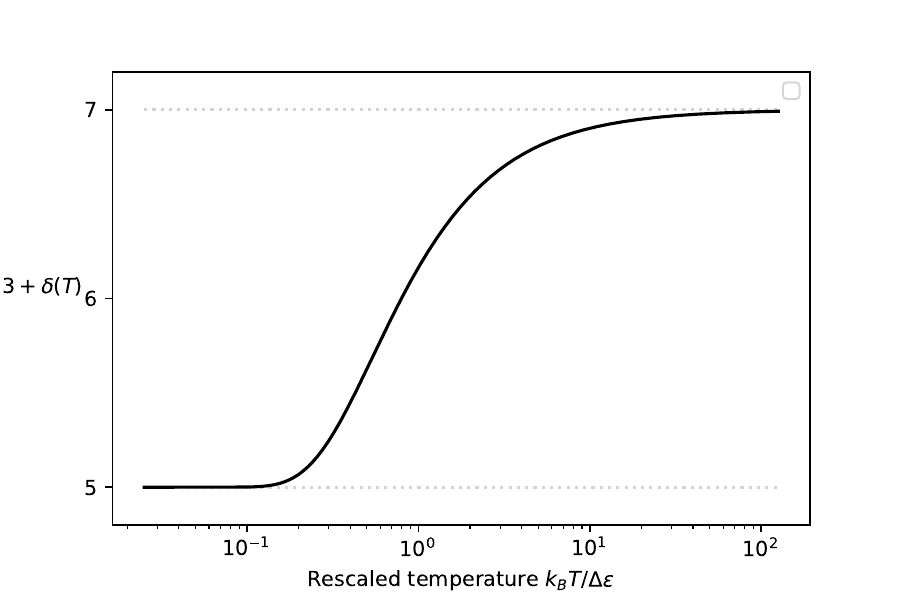}
    \caption{Plot of the total number of degrees of freedom $3 + \delta$ as a function of the rescaled temperature $k_B T / \Delta \epsilon$.}
    \label{figure:the_delta}
    \end{subfigure} \hspace{5pt}
\centering
\begin{subfigure}[t!]{0.48\linewidth}
\includegraphics[width=\linewidth]{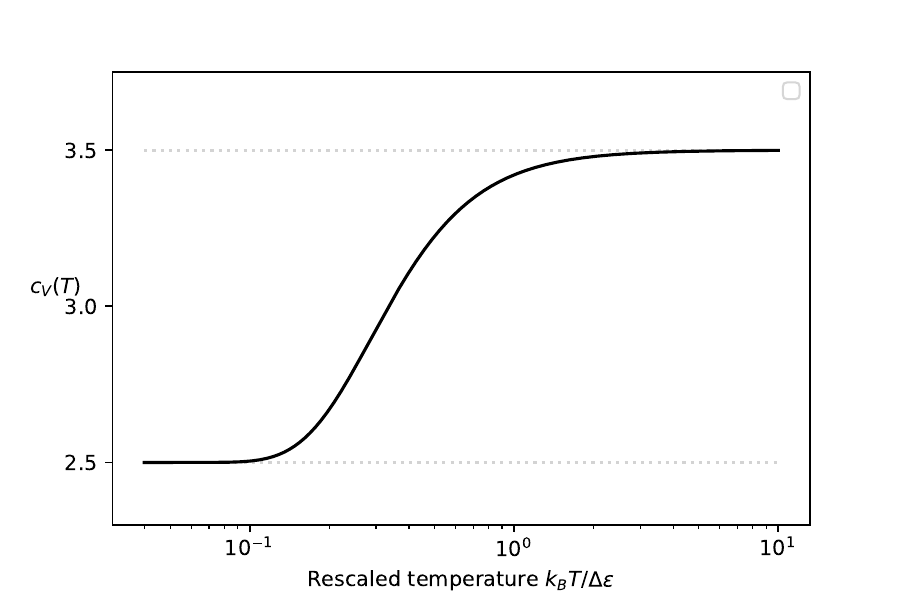}
    \caption{Plot of the heat capacity at constant volume $c_V$ as a function of the rescaled temperature $k_B T / \Delta \epsilon$.}
    \label{figure:the_cv}
\end{subfigure}
    \caption{Plot of the total number of degrees of freedom and heat capacity at constant volume associated with the model described in Sec. \ref{section2}.}
    \label{figure:d_and_c}
\end{figure}

\noindent We plot these functions, defined by~\eqref{eq:delta}--\eqref{eq:cv}, in Fig.~\ref{figure:d_and_c}.

\noindent {\bf Acknowledgements} 
Work performed in the frame of activities sponsored by University of Parma, CERMICS of École Nationale des Ponts et Chaussées, and Italian INdAM-GNFM. MB and MG also thank the support of the project PRIN 2022 PNRR
{\em Mathematical Modelling for a Sustainable Circular Economy in Ecosystems} (project code P2022PSMT7, CUP D53D23018960001) funded by the European Union - NextGenerationEU PNRR-M4C2-I 1.1 and by MUR-Italian Ministry of Universities and Research, and also of the action ``Bando di Ateneo 2022 per la ricerca'' co-funded by
MUR-Italian Ministry of Universities and Research - D.M. 737/2021 - PNR - PNRR - NextGenerationEU, Project ``Collective and Self-Organised Dynamics: Kinetic and Network Approaches''. 

\newpage

\bibliographystyle{plain}
\bibliography{biblio}

\end{document}